# Optimization models and algorithms for the Unit Commitment problem

**Javal Vyas[a], Carl Laird[a], Ignacio E. Grossmann[a,*], Ricardo M. Lima[b], Iiro Harjunkoski[c], Jan Poland[d]**

[a] Department of Chemical Engineering, Carnegie Mellon University, Pittsburgh, PA 15213, USA
[b] King Abdullah University of Science and Technology, Thuwal 23955-6900, Saudi Arabia
[c] Hitachi Energy Research, Havellandstr. 10-14, Mannheim 68309, Germany
[d] Hitachi Energy Research, Segelhofstr. 7, 5405 Baden-Dättwil, Switzerland
* Corresponding Author: grossmann@cmu.edu.

## ABSTRACT

The unit commitment problem determines the optimal strategy to meet the electricity demand at minimum cost by committing power generation units at each point of time. Solving the unit commitment problem gives rise to a challenging optimization problem due to its combinatorial complexity and potentially long solution time requirements. Our proposed solution approach utilizes a decomposition method in conjunction with alternative models from the EGRET library. Results of this decomposition approach tested against four benchmarking systems show that significant computational speed ups are achieved.

## BACKGROUND

The Unit Commitment (UC) problem is a foundational problem in the field of power systems. This problem tackles the decision making for the operation of the grid, defining when to switch on or off generators, and defining their power production levels to meet the forecasted energy demand. Initially it was solved using dynamic programming [1], which moved to a more sophisticated version using the heuristics based on Lagrangean Relaxation and then to mixed-integer linear programming (MILP) formulations[2]. This problem is usually solved for the day-ahead unit commitments, and bids are placed by the operator, the downside to this is that the operator has at most around 3 hours to solve the model and analyze the bids of which the most time is consumed in data verification, network analysis. Thus, the solution time requirement for solving the UC problem should be no more than 10-15 minutes. In practice, the UC problem must be solved repeatedly throughout the day, and the practical time window for the entire decision-making process is often around 20-30 minutes. This creates a stringent need for a faster solution approach. As the number of generation units increases due to the smaller distributed electricity resources (often solar and wind), this also increases the problem complexity significantly.

## PROBLEM STATEMENT

In the context of an electrical grid comprising multiple generation units, this study addresses the security constrained unit commitment (SCUC) problem, which also considers power grid transmission line limitations and other operational constraints such as reserve requirements. The objective is to determine the optimal schedule for generation units over a defined time horizon. Here we deal with the most common day-ahead unit commitment problem, and thus the time horizon remains 24 hours with hourly intervals. The grid's configuration includes the operational parameters of each generation unit, network topology, transmission line limits, and the forecasted energy demands at each node. The aim is to minimize the total generation cost, while satisfying to the energy balance constraints at each node, transmission line capacity limits, and operational limits of the generation units.

As the problem size is expected to rapidly grow due to both the increased number of generation units, as well as a trend to move to 15-minute time intervals, this work specifically focuses on efficient modeling and solution of a large-scale SCUC problem. It excludes considerations of uncertainties associated with renewable energy sources, as well as the inclusion of energy storage. Here it is assumed that the fuel sources for all generators are non-renewable and thus fully dispatchable.

## PROPOSED METHODOLOGY

Given the extensive body of research on UC problem, this study leverages established tools and methodologies from the literature. We employ the open-source software package EGRET[3], which provides a comprehensive framework for UCP by offering 41 different model formulations drawn from the existing literature.

Among these formulations, we focus on two specific categories of models: those developed by Morales-España et al. [4], Knueven et al. [5], known for their competitive computational performance, and the models identified by Knueven et al. [3], classified as the "Tight" and "Compact" formulations based on their efficiency in computational studies.

However, these models often struggle to close the optimality gap within a reasonable time limit, when applied to larger, more complex systems. To address this limitation, we propose the use of a shrinking horizon method.

### Optimization models

To expand on the available models from the literature, we briefly examine the key characteristics and improvements introduced by each model:

a. Morales-España model (MLR): Building upon the foundational model developed by Carrion and Arroyo[6], the Morales-España model introduces tighter bounds on generation production. The proposed reformulation of the MILP model is both tight and compact, thus making it computationally competitive.

b. Knueven Ostrowski Watson model (KOW): The KOW model refines the Morales-España model by tightening the formulations of both piecewise linear production costs and start-up costs. The proposed model aggregates the thermal generators with identical characteristics into a single unit which reduces the redundancy in search space.

c. Tight Model (T): The tight model, akin to the KOW model improves on the tightness by enhancing the generation limits, tight start-up cost formulation, and it also uses two period ramping from the Damci Kurt et al. [7]. The tight model does lose on the compactness. Here the tightness of the model refers to the search space a solver needs to explore to find the optimal solution. Tighter models would have reduced search space.

d. Compact model (Co): The Compact model, is similar to the MLR model, with the difference being two period ramping from Damci Kurt et. al. [7] and improved reserve production. This formulation loses a bit on the tightness but is compact. Here the compactness of the model refers to the model size. A compact model has fewer number of variables and constraints.

### Shrinking Horizon Method

The Shrinking Horizon method was introduced by Balasubramanian and Grossmann[8], is an iterative method designed to improve computational tractability of multiperiod MILP models, particularly for large-scale systems. In this approach, the time horizon is divided into overlapping windows. During the first iteration, a window spanning from time t=1 to t=h is defined, where binary decision variables are active. For the remaining time-periods t=h+1 to t=n, the binary variables are relaxed. The model is solved within this window, and the binary decisions are recorded.

In subsequent iterations, the binary decisions for the previous window are fixed, and the window is advanced to t=h+1 to t=2h. This process is repeated iteratively, shifting the window forward until the entire time horizon is covered. This method allows for a more manageable and efficient solution process as illustrated in Figure 1.

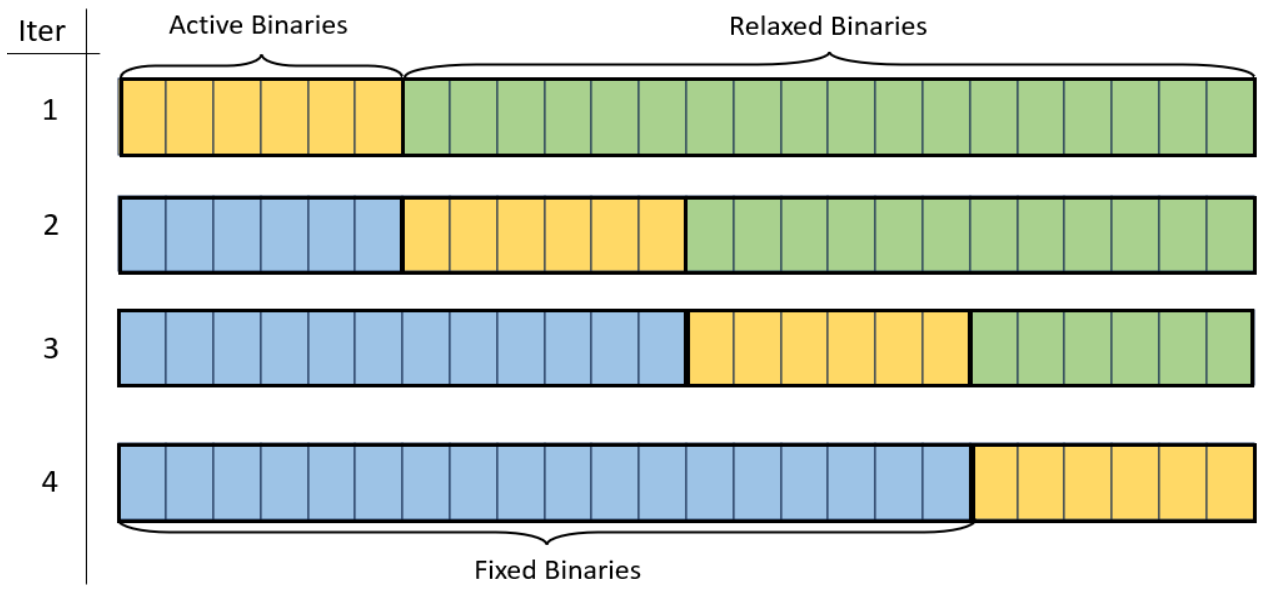

**Figure 1.** Shrinking Horizon Method

The Shrinking Horizon method inherently relies on solving a sequence of optimization problems over progressively shorter horizons, where future time periods are relaxed to facilitate computational tractability. However, the effectiveness of this approach depends on the quality of these relaxations, loose relaxations can lead to poor bounds, increased solution time, and suboptimal intermediate decisions.

Tight Model, mentioned to have the tightest relaxation amongst formulations in the UC literature [3] is particularly well-suited for the Shrinking Horizon method. The tight relaxations offer stronger bounds at each iteration and thus enhances the pruning efficiency of branch-and-cut algorithm, reducing the number of nodes explored in the tree search. Thus, the solvers can converge faster to a high-quality solution.

Furthermore, the Tight model mitigates the myopic effects of the Shrinking Horizon method by ensuring the relaxed future decisions remain closer to their integer counterparts. This reduces the risk of infeasibilities or drastic corrections in subsequent iterations, leading to a stable and efficient optimization process. By leveraging the Tight Model within the Shrinking Horizon framework, the full computational potential of modern off-the-shelf solvers can be effectively realized.

## CASE STUDY

In this section, we present the results of applying the Shrinking Horizon method across various benchmark case studies, including IEEE118, WP2383, SP3120, and CASE6468RTE, each representing a system of varying size and complexity (Table 1). These case studies range from IEEE118, a relatively small system with 17 generating units, to the largest system, CASE6468RTE, containing 1181 units. 20 instances for each case study using a proprietary code by our collaborators at Hitachi Energy[9]. The analysis focuses on comparing the performance of the Shrinking Horizon method with the traditional monolithic approach, evaluating both computational speedup and solution quality in terms of objective deviation.

All case studies in the monolithic approach were solver to 0.01% optimality gap, except for CASE6468RTE, which was solved to a 1% gap and time limit of 3600

**Table 1:** Network Topology

| Case name | Nodes (Buses) | Lines | Generating Units |
|---|---|---|---|
| IEEE118 | 118 | 179 | 17 |
| WP2383 | 2383 | 2896 | 293 |
| SP3120 | 3120 | 3693 | 452 |
| CASE6468RTE | 6468 | 9000 | 1181 |

**Table 2:** Case6468RTE Analysis

| | Compact Model | | KOW Model | | MLR Model | | Tight Model | |
|---|---|---|---|---|---|---|---|---|
| | Optimalilty Gap % | Speedup factor | Optimalilty Gap % | Speedup factor | Optimalilty Gap % | Speedup factor | Optimalilty Gap % | Speedup factor |
| Min | **0.0000** | **0.5312** | **0.0000** | 0.4834 | **0.0000** | 0.3776 | **0.0000** | 0.4386 |
| Median | 4.2662 | 0.9455 | 6.0975 | 1.7521 | 5.0206 | 0.9491 | **4.1889** | **1.9803** |
| Max | 80.5953 | 1.8351 | 75.7198 | **4.2460** | 71.9478 | 1.9836 | **71.1481** | 2.8829 |
| Std | 17.8190 | 0.3989 | 24.6893 | 0.8674 | 20.6174 | 0.4739 | 20.1230 | 0.6995 |

**Table 3:** SP3120 Analysis

| | Compact Model | | KOW Model | | MLR Model | | Tight Model | |
|---|---|---|---|---|---|---|---|---|
| | Optimalilty Gap % | Speedup factor | Optimalilty Gap % | Speedup factor | Optimalilty Gap % | Speedup factor | Optimalilty Gap % | Speedup factor |
| Min | 0.0006 | 0.4002 | **0.0001** | 0.5539 | 0.0012 | 0.4611 | 0.0149 | **0.6022** |
| Median | **0.0766** | 0.6839 | 0.1967 | 1.2206 | 0.1069 | 0.7399 | 0.1785 | **1.3568** |
| Max | 100.0 | 2.2878 | 0.7198 | 3.9246 | **0.5909** | 2.9072 | 0.6640 | **4.0035** |
| Std | 22.325086 | 0.5263 | 0.1934 | 0.9948 | 0.1961 | 0.6610 | 0.1997 | 1.0178 |

**Table 4:** WP2383 Analysis

| | Compact Model | | KOW Model | | MLR Model | | Tight Model | |
|---|---|---|---|---|---|---|---|---|
| | Optimalilty Gap % | Speedup factor | Optimalilty Gap % | Speedup factor | Optimalilty Gap % | Speedup factor | Optimalilty Gap % | Speedup factor |
| Min | 0.0021 | 0.4455 | 0.0025 | **0.5468** | 0.0021 | 0.4611 | **0.0012** | 0.4653 |
| Median | **0.0594** | 0.6891 | 0.0651 | **1.1884** | 0.0621 | 0.7399 | 0.0792 | 0.8688 |
| Max | 0.2087 | 2.0498 | 0.2433 | **3.6665** | **0.1266** | 2.9072 | 0.2073 | 3.4818 |
| Std | 0.0592 | 0.4092 | 0.067 | 0.8085 | 0.0451 | 0.6443 | 0.0593 | 0.7771 |

**Table 5:** IEEE118 Analysis

| | Compact Model | | KOW Model | | MLR Model | | Tight Model | |
|---|---|---|---|---|---|---|---|---|
| | Optimalilty Gap % | Speedup factor | Optimalilty Gap % | Speedup factor | Optimalilty Gap % | Speedup factor | Optimalilty Gap % | Speedup factor |
| Min | 0.0000 | 0.4044 | **0.0000** | 0.5568 | **0.0000** | **0.6685** | 0.0000 | 0.5058 |
| Median | 0.0002 | 1.0759 | **0** | **2.3391** | **0** | 1.8660 | 0.0185 | 1.4557 |
| Max | 5.7381 | 2.3007 | **1.6808** | **11.0642** | **1.6808** | 3.7002 | 4.3298 | 5.5728 |
| Std | 1.5072 | 0.5464 | 0.3894 | 2.2481 | 0.4074 | 0.9062 | 1.3254 | 1.2035 |

* **Best** values are bold and underlined. Second best values are underlined

seconds. For most part of the CASE6468RTE case study in the monolithic approach, the solver timed out with a median optimality gap of nearly 5%. Notably, the Shrinking Horizon method was solved for 0% optimality gap for each iteration across all case studies, except for CASE6468RTE, where each iteration was solved to 1% optimality gap with a time limit of 3600 seconds per iteration.

The computational study was conducted on a Linux machine equipped with 8 Intel Xeon Gold 6234 CPUs, operating at 3.30 GHz with a total of 8 hardware threads and 1 TB of RAM, all running within an Ubuntu environment. The models were constructed using EGRET's model library and solved with Gurobi v9.5.1 as the solver.

The speedup factor is calculated as:

$$speedup\ factor = \frac{solution\ time\ with\ monolithic\ model}{solution\ time\ with\ SH\ method}$$

The speedup factor reflects the efficiency gained by applying the shrinking horizon method. A factor greater than 1 suggests that the SH method is less effective than solving the problem as a whole.

Similarly, the objective deviation factor is calculated as:

$$Obj\ deviation\ factor\ \% = \frac{Obj\ value\ w/SH}{Obj\ value\ w/o\ SH}$$

This factor behaves conversely to the speedup factor: a higher objective deviation factor indicates a greater deviation from the optimal solution, meaning the SH method is less effective and results in a compromise on the objective value.

Similarly, a higher optimality gap signifies a larger deviation from the optimal solution, implying a greater trade-off in solution quality for the computational efficiency achieved with the shrinking horizon method.

The results demonstrate that the solving times for the shrinking horizon method are consistently lower, and in some cases, significantly lower, for the KOW and Tight models across most case studies (Figure 2). This indicates that these models are more effective under the shrinking horizon method compared to other models. As a result, they emerge as the preferred formulations for achieving the best speedups.

**Figure 2.** Casewise analysis of models

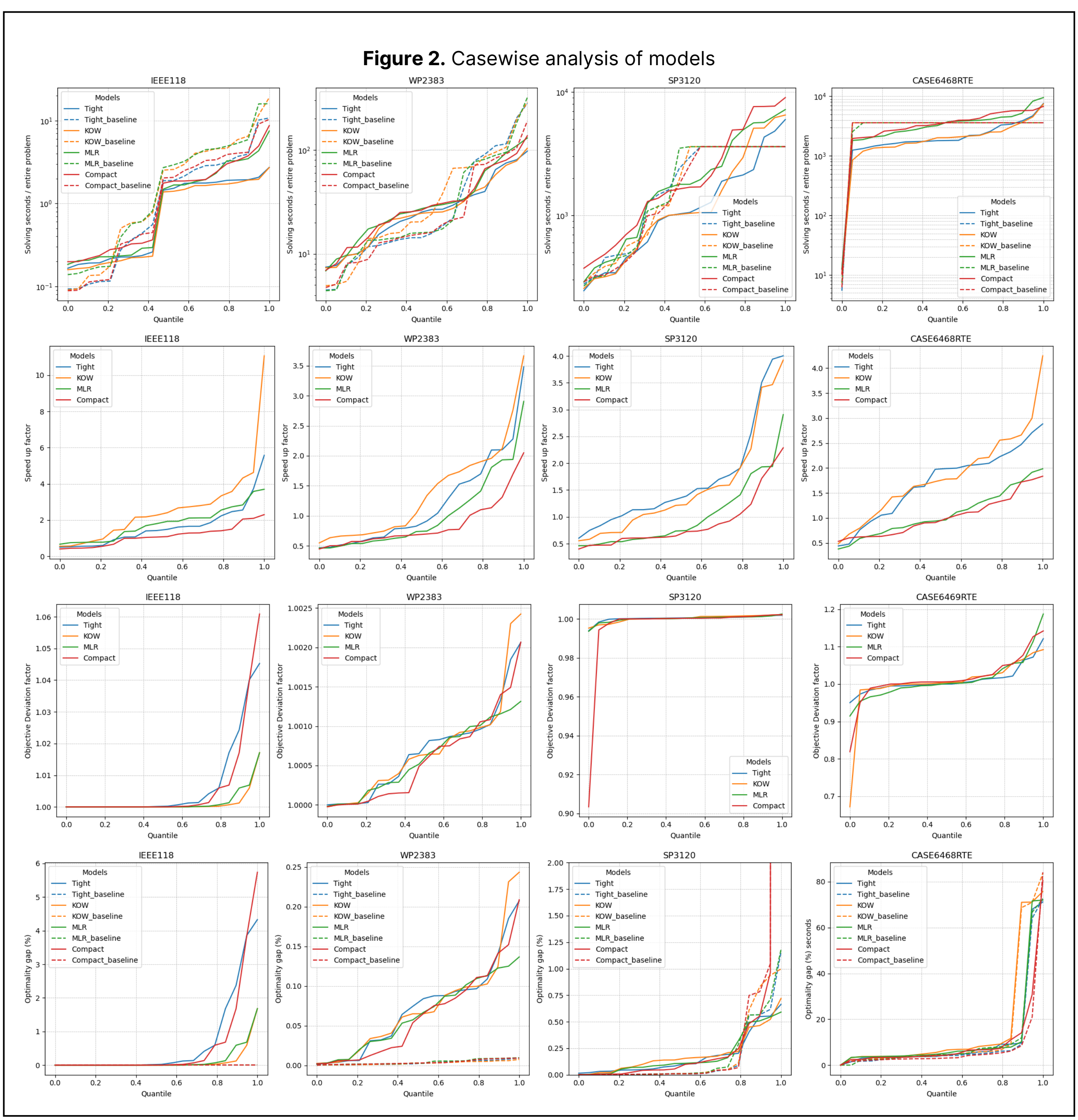

Additionally, these models maintain near-optimal performance, with minimal compromise on the objective value. For example, the maximum median deviation in the objective value is approximately 0.25% for both the Tight and KOW models for all other cases except case6468rte where the median deviation is 4.18% for Tight model while KOW has 4.26% as the median deviation, demonstrating their ability to balance computational efficiency with solution accuracy.

The advantage of these models becomes even more pronounced when applied to large case studies, where their efficiency and reliability provide substantial benefits. This highlights their suitability for scaling to more complex or computationally demanding scenarios.

## CONCLUSION

In this study, we have proposed and validated a novel approach to solving the UC problem by integrating the Tight, KOW, MLR and compact model with a Shrinking Horizon method.

Our analysis of 20 random instances for each case study further confirmed the Tight and KOW model's stability and computational efficiency, particularly excelling in larger systems like CASE6468RTE. These results demonstrate that the Tight and KOW model with the Shrinking Horizon method are effective strategies for efficiently solving large-scale UCPs while ensuring computational speedups while not losing much on the solution quality.

Future research should explore this approach's scalability in more complex grid scenarios involving renewable energy sources, storage systems, and uncertainties. Additionally, investigating alternative decomposition techniques could further enhance its performance and applicability.

## ACKNOWLEDGEMENTS

The authors would like to acknowledge the financial support for this research from Hitachi Energy Research and the Center for Advanced Process Decision-making at Carnegie Mellon University.

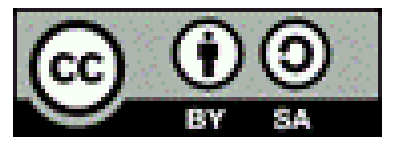